\documentclass[12pt]{amsart}

\usepackage{amsmath, amsthm, amssymb}
\usepackage[all,cmtip]{xy}
\usepackage{hyperref}

\newcommand{\End}{\text{End}}
\newcommand{\pdeg}{\ensuremath{\text{pdeg}}}
\newcommand{\rk}{\ensuremath{\text{rk}}\, }
\newcommand{\id}{\ensuremath{\text{Id}}}
\newcommand{\parend}{\ensuremath{\mathcal{P}ar\mathcal{E}nd}}
\newcommand{\sparend}{\ensuremath{\mathcal{SP}ar\mathcal{E}nd}}
\newcommand{\fparend}{\ensuremath{\mathcal{FP}ar\mathcal{E}nd}}
\newcommand{\phig}{\ensuremath{\mathcal{P}_{\alpha}}}
\newcommand{\pvect}{\ensuremath{\mathcal{N}_{\alpha}}}

\newcommand{\pfram}{\ensuremath{\mathcal{F}_{\alpha}}}
\newcommand{\pmast}{\ensuremath{\mathcal{F}}}
\newcommand{\im}{\ensuremath{\text{im}}}

\newcommand{\git}{  / \! \! /}
\newcommand{\Gr}{\ensuremath{\text{Gr}}}

\newtheorem{theorem}{Theorem}[section]
\newtheorem{lemma}[theorem]{Lemma}
\newtheorem{proposition}[theorem]{Proposition}

\newtheorem{definition}[theorem]{Definition\rm}

\DeclareMathOperator{\GL}{GL\,}

\begin{document}

\title{Moduli of Parabolic Higgs Bundles and Atiyah Algebroids}
\author{Marina Logares}
\address{Departamento de Matem\'atica Pura, Faculdade de Ci\^encias, Universidade do Porto, Rua do Campo Alegre, 687, 4169-007 Porto, Portugal}
\curraddr{Departamento de Matem\'aticas, CSIC, Serrano 121, 28006 Madrid, Spain}
   \email{marina.logares@icmat.es}
\author{Johan Martens}
\address{Department of Mathematics, University of Toronto, 40 St. George Street, Toronto Ontario M5S 2E4, Canada}
\curraddr{Centre for Quantum Geometry of Moduli Spaces,
Department of Mathematical Sciences, Aarhus University,
Ny Munkegade 118, bldg. 1530,
DK-8000 \r{A}rhus C,
Denmark }
   \email{jmartens@imf.au.dk}
\date{\today}
\setcounter{tocdepth}{3}

\begin{abstract}In this paper we study the geometry of the moduli space of (non-strongly) parabolic Higgs bundles over a Riemann surface with marked points.  
We show that this space possesses a Poisson structure, extending the one on the dual of an Atiyah algebroid over the moduli space of parabolic vector bundles.
By considering the case of full flags, we get a Grothendieck-Springer resolution for all other flag types, in particular for the moduli spaces of twisted Higgs bundles, as studied by Markman and Bottacin and used in the recent work of Laumon-Ng\^o.  We discuss the Hitchin system, and demonstrate that all these moduli spaces are integrable systems in the Poisson sense.
\end{abstract}

\maketitle

\tableofcontents

\section{Introduction}
Higgs bundles, introduced by Hitchin in \cite{hitchin1, hitchin2}, have e\-mer\-ged in the last two decades as a central object of study in geometry, with several links to physics and number theory.  Over a smooth compact Riemann surface the moduli space of Higgs bundles contains as a dense open subset the total space of the cotangent bundle to the moduli space of vector bundles.  In fact the induced complex symplectic form is part of a hyper-K\"ahler structure and extends to the whole of the moduli space of Higgs bundles, and it is a celebrated fact that the moduli space comes equipped with a algebraically completely integrable system, through the Hitchin map.  
\\

A natural generalization of vector bundles arises when one endows the vector bundle with a parabolic structure~\cite{MS}, i.e. with choices of flags in the fibers over certain marked points on the Riemann surface.  One can talk of Higgs bundles in that setting as well, as was first done by Simpson in~\cite{simpson}.  Various choices can be made for this. In order to have the corresponding moduli space contain as an open subset the total space of the cotangent bundle to the moduli space of parabolic vector bundles, replicating the non-parabolic situation, several authors~\cite{michael,konno,bettipar, faltings} restrict the parabolic Higgs bundles to those that we shall refer to as strongly parabolic, meaning that the Higgs field is nilpotent with respect to the flag.  
\\

One can however also demand the Higgs field to simply respect the parabolic structure at the marked points, and a moduli space $\phig$ for those was constructed by Yokogawa in~\cite{yoko2}.  The locus of that moduli space where the underlying parabolic vector bundle is stable again forms a vector bundle over the moduli space of parabolic vector bundles $\pvect$.  We show here that this vector bundle is the dual of an Atiyah algebroid associated with a principal bundle over $\pvect$, the structure group for which is the product of the Levi groups given by the various flags at the marked points, modulo $\mathbb{C}^*$ to account for global endomorphisms of the bundle. As the dual of an algebroid its total space carries a complex algebraic Poisson structure, which in fact extends to the whole of $\phig$.  In the particular situation where all the flags are trivial, $\{0\}\subset E|_p$, this was already shown independently by Bottacin~\cite{bottac} and Markman~\cite{markman}.  We further study the Hitchin system for $\phig$ and its symplectic leaves, showing that this makes the $\phig$ for all flag types into integrable systems in the sense of Poisson geometry, with Casimir functions that generically induce the foliation of symplectic leaves, which are integrable systems in the usual symplectic sense.
Though we don't explicitly use Lie groupoids, our philosophy is very much that the symplectic leaves are the co-adjoint orbits for the groupoid determined by the principal bundle.
\\

With this in mind we also look at the forgetful morphisms between such moduli spaces of various flag types.  We show that they are Poisson and generically finite.  By looking at such morphisms starting from the moduli space for full flags we obtain a global analogue of the Gro\-then\-dieck-Springer resolution of Lie algebras, as the moduli space for full flags is a regular Poisson variety.  As the Grothendieck-Springer resolution plays a crucial role in modern geometric representation theory, this opens perspectives on generalizing classical constructions to this global setting\footnote{Indeed, we very recently became aware of the preprint~\cite{compet}, where parabolic Higgs bundles are used to generalize the Springer theory of Weyl group representations to a `global' setting, without however taking our viewpoint of Atiyah algebroids and groupoids.}.
\\

The presentation we have given is largely done in the language of Lie algebroids, but one could have reformulated everything we say about Atiyah algebroids in terms of Poisson reduction of cotangent bundles.  The choice is partly one of personal preference, and partly due to the fact that the Atiyah sequence of the algebroid naturally follows from the deformation theory of parabolic vector bundles.  Similarly, though almost all of our arguments use hypercohomology, we have avoided the use of derived categories, in order to exhibit the Poisson structures more explicitly.  
\\

We remark here that our entire construction is depending on the existence of the principal bundle over $\pvect$.  In the appendix~\ref{append} we describe a construction of this bundle for full flags, using previous work by Hurtubise, Jeffrey and Sjamaar~\cite{framparshea}, and outline a possible construction for other flag types.
\\

This paper is organized as follows: in section~\ref{intro} we give the necessary background regarding parabolic Higgs bundles and their moduli, as well as a description of the Hitchin fibration. Most of this is standard, with the possible exception of the observation (Proposition~\ref{unique}) that the smoothness of the spectral curve implies that the Higgs field uniquely determines the parabolic structure, even when eigenvalues are repeated.  In section~\ref{poissonstruc} we give the required background material regarding Lie groupoids and Lie algebroids, and prove the main result of this article, the interpretation of the moduli space as a partial compactification of the dual of an Atiyah algebroid.  We also show here that for the induced Poisson structure the Hitchin map is an integrable system.  In section~\ref{G-S} we remark that for nearby parabolic weights the morphisms between the various moduli spaces (with different flag structures) are Poisson, giving a Grothendieck-Springer resolution by means of the full flags. In section~\ref{further} we discuss the relationship of our work with the earlier results by Bottacin and Markman, as well as further directions.  Appendix~\ref{append} discusses a construction of the principal bundle over the moduli space of parabolic vector bundles used in the main theorem.

\subsection{Remark on notation}  Unfortunately nomenclature conventions regarding parabolic Higgs bundles vary in the literature.  For us a parabolic Higgs bundle will only require the Higgs field at a marked point to respect the flag there.  We will refer to the special case where the Higgs field is nilpotent with respect to the filtration as a \emph{strongly} parabolic Higgs bundle.

\subsection{Acknowledgements} The authors would like to thank Sergey Ar\-khi\-pov, David Ben-Zvi, Philip Boalch, Hans Boden, Chris Brav, Ron Donagi, Tom\'as Gom\'ez, Peter Gothen,  Tam\'as Hausel, Nigel Hitchin,  Jacques Hurtubise, Lisa Jeffrey, Raj Mehta, Eckhard Meinrenken, Szi\-l\'ard Szab\'o and Michael Thaddeus for useful conversations, remarks, and encouragement, as well as the MPI Bonn, NSERC and the Centro de Matem\'atica da Universidade do Porto for financial support.

\section{Moduli spaces of Parabolic Higgs bundles} \label{intro} 
\subsection{Parabolic vector bundles}
Let $X$ be a compact Riemann surface or smooth complex\footnote{Presumably all the results below hold over arbitrary algebraically closed fields.  Neither Yokogawa's construction of the moduli spaces we shall use, nor any of our work, requires the ground field to be $\mathbb{C}$.  We do rely crucially however on the results of \cite{BNR}, which assumes the characteristic to be zero.} projective curve of genus $g$ with $n$ distinct marked points $p_1,\dots, p_n$. If $g=0$ we assume $n\geq 3$, if $g=1$ we assume $n\geq 1$. Let $D$ be the effective reduced divisor $p_1+\ldots+p_n$.  A \emph{parabolic vector bundle} (\cite{MS}) on $X$ is an algebraic rank $r$ vector bundle $E$ over $X$ together with a parabolic structure, i.e. a (not necessarily full) flag for the fiber of $E$ over the marked points
$$E|_{p}=E_{p,1}\supset \cdots \supset E_{p,r(p)}\supset \{0\},$$
together with a set of parabolic weights
$$0\leq \alpha_{1}(p)<\cdots < \alpha_{r(p)}(p)<1.$$
We denote the multiplicities by $m_{i}(p)=\dim E_{p,i} - \dim E_{p, i+1}$, and the associated graded as $\Gr (p)=\oplus_i E_{p,i}/E_{p,i+1} $.  

Note that the structure group of the bundle $E$ is  $\GL(r)$.  In terms of the associated frame bundle the parabolic structure corresponds to a reduction of the structure group of this principal bundle to a certain parabolic subgroup of $\GL(r)$ at each marked point $p_i$. 
 We will denote this parabolic subgroup by $P_p$, and its corresponding Levi group by $L_p$, with Lie algebra $\mathfrak{l}_p$. For the sake of convention we will fix as a Borel subgroup in $GL(r)$ the lower triangular matrices, and all parabolic subgroups are taken to contain this Borel.
 \\
 
 We will further need linear endomorphisms\footnote{A priori we don't require morphisms between vector bundles to have constant rank.}  $\Phi$ of a parabolic vector bundle which are either parabolic - meaning that at the fiber over a marked point $p$ we have $\Phi|_p(E_{p,i}) \subset E_{p,i}$ - or \emph{strongly parabolic} - meaning that $\Phi|_p(E_{p,i}) \subset E_{p,i+1}$.  We denote the sheaves of parabolic respectively strongly parabolic endomorphisms as $\parend(E)$ and $\sparend(E)$.
\\

The relevance of the parabolic weights $\alpha$ comes  from  the notion of \emph{parabolic degree} of a bundle, denoted by \emph{$\pdeg$} :
$$\pdeg(E) = \deg(E) + \sum_{p\in D} \sum_i m_i(p) \alpha_{i}(p),$$
which satisfies the Gauss-Chern formula for connections with logarithmic singularities (see Proposition 2.9 in \cite{biquard}).  The $\alpha$ also occur in the celebrated Mehta-Seshadri theorem (\cite{MS}) that establishes a correspondence between stable parabolic bundles and unitary representations of the fundamental group of the punctured surface $X\setminus D$, where they determine the holonomy around the punctures.
\\

Every algebraic subbundle $F$ of $E$ is naturally given the structure of a parabolic bundle as well, by simply intersecting $F|_p$ with the elements of the flag of $E|_p$,  discarding any such subspace of $F|_p$ that coincides with a previous one, and endowing it with the largest of the corresponding parabolic weights:
$$\textrm{for}\, F_{p,i}=F|_p\cap E_{p,j}\qquad \alpha^{F}_i(p)=\max_{j}\{\alpha_j ;\ F|_p\cap E_{p,j}=F_{p,i}\}.$$

We say that a parabolic vector bundle is stable  if for each proper subbundle $F$ we have that \begin{equation}\label{slopeineq}\frac{\pdeg(F)}{\rk(F)}< \frac{\pdeg(E)}{\rk(E)}.\end{equation}  Semi-stability is defined similarly, by asking for weak inequality.  
The weights are called \emph{generic} when stability and semi-stability coincide. Note that the term generic is used in the sense that the set of non-generic weights has positive codimension.
\\

There exists a moduli space for semi-stable parabolic vector bundles~\cite{MS}, which we shall denote by $\pvect$.  This is a normal projective variety of dimension $$\dim \pvect = (g-1)r^2+1+ \sum_{p\in D}\frac{1}{2}\left(r^2-\sum_i m_{i}(p)^2\right)$$ and when the weights are generic it is non-singular.  From now on we will assume genericity of weights, even though in the non-generic case all of what we say can still be carried through when restricted to the stable locus of $\pvect$.
\\

We will need a further lemma that states that stability implies simplicity.  This is completely standard in the non-parabolic setting (see e.g. \cite[Corollary to Proposition 4.3]{narasesh}) but a parabolic version does not seem to have appeared in the literature, so we include it here for completeness.

\begin{lemma}\label{stablesimple}
If $E$ is an $\alpha$-stable parabolic vector bundle over $X$, then $$H^0(X, \parend(E))=\mathbb{C}$$ and $$H^0(X, \sparend(E))=\{0\}.$$
\end{lemma}
\begin{proof} 
First note that given a rank $r$ algebraic vector bundle over a complete non-singular variety, the $r$-th power of any endomorphism of this bundle necessarily has constant rank, since the coefficients of the characteristic polynomial are regular functions and hence constant and the $r$-th power has no nilpotent Jordan blocks.  Therefore we can talk of the kernel of this $r$-th power of the endomorphism as algebraic subbundles.

Now, given a parabolic endomorphism $e$ of a parabolic vector bundle $E$, $e$ is zero or an isomorphism if and only if $e^r$ is.  We shall therefore consider $f=e^r$.  The subbundle $\ker(f)$ has a canonical induced parabolic
structure.
The same is true for $\im(f)$, which we also think of as a subbundle of $E$.  The parabolic weights that $\ker(f)$ and $\im(f)$ inherit as subbundles of $E$, when counted with multiplicities, are complementary to each other with respect to the parabolic weights of $E$.

Assume now that $f$ is neither zero nor an isomorphism, so both $\ker(f)$ and $\im(f)$ are proper subbundles of $E$.  Using the stability we have that 
\begin{equation}\label{mar1}
\frac{\pdeg\left(\ker(f)\right)}{\rk\left(\ker(f)\right)}
<
\frac{\pdeg(E)}{\rk(E)}
\text{ and }
\frac{\pdeg\left(\im(f)\right)}{\rk\left(\im(f)\right)}
<
\frac{\pdeg(E)}{\rk(E)}.
\end{equation} 

Using the complementarity of the parabolic weights of $\ker(f)$ and $\im(f)$ however we also have that
\begin{equation}\label{mar2}
\frac{\pdeg(E)}{\rk(E)}=
\frac
{\pdeg\left(\ker(f)\right)+\pdeg\left(\im(f)\right)}
{\rk\left(\ker(f)\right)+\rk\left(\im(f)\right)}.
\end{equation}
One easily sees that the combination of (\ref{mar1}) and (\ref{mar2}) would give $$\frac{\pdeg\left(\ker(f)\right)}{\rk\left(\ker(f)\right)}<\frac{\pdeg\left(\im(f)\right)}{\rk\left(\im(f)\right)} \text{ and } \frac{\pdeg\left(\im(f)\right)}{\rk\left(\im(f)\right)} < \frac{\pdeg\left(\ker(f)\right)}{\rk\left(\ker(f)\right)}, $$
hence $f$ is either zero or an isomorphism.

If $f$ is an isomorphism, one just has to take a point $x\in X$ and consider an
eigenvalue $\lambda$ of $f_x:E_x\to E_x$. Next consider the parabolic endomorphism
of $E$ given by  $(f - \lambda\, \id_E)$, by the same reasoning as before this is an
isomorphism
or it is zero. Hence $H^0(\parend(E))=\mathbb{C}$.

In the case of a strongly parabolic endomorphism, as for any marked point $p\in D$ zero is
one of the eigenvalues of $f_p$, one gets similarly\\ $H^0(\sparend(E))=\{ 0\}.$\end{proof}

Finally, notice that $\parend{E}$ is naturally dual to $\sparend(E)(D)$, and vice versa $\sparend(E)$ is dual to $\parend(E)(D)$.  Throughout the paper we shall often use Serre duality for the hypercohomology of a complex on a curve, so we recall its statement: for a bounded complex $\mathcal{C}$ of locally free sheaves on $X$ of the form $$0\rightarrow C^0\rightarrow \ldots \rightarrow C^m \rightarrow 0$$ we have the natural duality $$\mathbb{H}^i(\mathcal{C})^*\cong \mathbb{H}^{1-i+m}(\mathcal{C}^*\otimes K),$$
where $\mathcal{C}^* \otimes K$ is the complex $$0\rightarrow (C^m)^*\otimes K \rightarrow \ldots \rightarrow  (C^0)^*\otimes K \rightarrow 0,$$
 again considered in degrees $0$ through $m$.

\subsection{Parabolic Higgs bundles}A \emph{parabolic Higgs bundle} (\cite{simpson}) is a parabolic vector bundle together with a \emph{Higgs field} $\Phi$, a bundle morphism $$\Phi: E \rightarrow E\otimes K(D) ,$$ where $K$ is the canonical bundle of $X$, which preserves the parabolic structure at each marked point: $$\Phi |_{p}(E_{p,i})\subset E_{p,i}\otimes K(D)|_{p},$$ i.e. $\Phi \in H^0(X, \parend(E)\otimes K(D))$.  
 In keeping with the notation introduced above we refer to the Higgs bundle as \emph{strongly} parabolic if the Higgs field is actually nilpotent with respect to the filtration, i.e. if
 $$\Phi|_{p}(E_{p,i})\subset E_{p,i+1}\otimes K(D)|_{p}.$$

Similar to vector bundles a Higgs bundle is (semi) stable if the slope condition $$\frac{\pdeg(F)}{\rk(F)}\underset{(=)}{<} \frac{\pdeg(E)}{\rk(E)}$$ holds, restricted now to all proper subbundles $F$ preserved by the Higgs field, i.e. with $\Phi(F)\subset F\otimes K(D)$.
\\

Denote by $\phig$ the moduli space of $\alpha$- semi-stable parabolic Higgs bundles of degree $d$ and rank $r$, which was constructed by Yokogawa in~\cite{yoko2} and further discussed in~\cite{bodenyoko}.  This space is a normal, quasi-projective variety of dimension \begin{equation}\label{dimyoko}\dim \phig = (2g-2+n)r^2 +1.\end{equation} Observe that this dimension is independent of the flag-type at the marked points, in contrast to the dimension of $\pvect$.  Indeed, there is a natural partial orderings on the flag types, and if ${\mathcal{P}}_{\tilde{\alpha}}$ is the corresponding moduli space for a finer flag type and the weights $\tilde{\alpha}$ and $\alpha$ are close enough such that semi-stability is preserved by the forgetful functor, then there is a forgetful morphism \begin{equation}\label{laatste}{\mathcal{P}}_{\tilde{\alpha}}\rightarrow\phig\end{equation} which is generically finite - see also proposition~\ref{unique} below.

\subsection{The Hitchin fibration}
\subsubsection{Hitchin map}Just as for ordinary Higgs bundles \cite{hitchin1, hitchin2}, the parabolic Higgs bundles form an integrable system by means of the \emph{Hitchin map}, defined as follows.  Given a vector bundle $E$, any invariant, homogeneous degree $i$ polynomial naturally defines a map $$H^0(\End (E)\otimes K(D) )\rightarrow H^0(K(D)^i).$$  Now, take the elementary symmetric polynomials as a homogeneous basis of polynomials on $\mathfrak{gl}(r)$ invariant under the adjoint action of $GL(r)$, then the corresponding maps $a_i$ combine to give the \emph{Hitchin map}
$$h_{\alpha}:\phig\rightarrow \mathcal{H}$$ where the vector space $\mathcal{H}$ is the Hitchin space $$\mathcal{H}=H^0(X,K(D))\oplus H^0(X,K(D)^2)\oplus \dots \oplus H^0(X,K(D)^r).$$ The components of $h_{\alpha}$ are defined as follows: for any parabolic Higgs bundle $(E,\Phi)$ and for any $x\in X$, let $k \in K(D)|_x$.  Then we have that $$\det(k.\id_{E|_x} -\Phi|_x) = k^r + a_1(\Phi)(x) k^{r-1}+\dots +a_{r-1}(\Phi)(x) k  + a_r(\Phi)(x),$$  and $h_{\alpha}(E,\Phi)$ is given by $(a_1(\Phi), \dots, a_r(\Phi) )$.  In~\cite{yoko2}, \S5,  it is shown that $h_{\alpha}$ is proper, and in fact projective.

Notice that $h_{\alpha}$ is \emph{blind to the parabolic structure} at each marked point, as it only depends on  $\Phi$ and the line bundle $K(D)$.  Indeed, suppose a given vector bundle $E$ and Higgs field $\Phi:E\rightarrow E\otimes K(D)$ can be equipped with two distinct parabolic structures compatible with $\Phi$, so as to obtain different stable parabolic Higgs bundles, possibly but not necessarily of different flagtype.  If we denote these two parabolic vector bundles by $F$ and $\widetilde{F}$, with $(F,\Phi)\in\mathcal{P}_{\alpha}$ and $(\widetilde{F},\Phi)\in \mathcal{P}_{\tilde{\alpha}}$, then necessarily $h_{\alpha}(F,\Phi)=h_{\tilde{\alpha}}(\widetilde{F},\Phi)$.

\subsubsection{Spectral curves}
For each element $s=(s_1,\dots, s_r)$ of $\mathcal{H}$ one can define a spectral curve $X_s$ in $S$, the total space of $K(D)$, as follows:  pull back $K(D)$ to $S$ and denote its canonical section as $\lambda$.  Then $X_s$ is the zero-locus of $$\lambda^r+s_1\lambda^{r-1}+\dots+s_r,$$ a (possibly ramified) covering of $X$.  As usual, by a Bertini argument, for a generic element in $\mathcal{H}$ the corresponding $X_s$ is smooth.  The genus of $X_s$ can be given using the adjunction formula: \begin{align*}2g(X_s)-2 &= \deg(K_{X_s})\\ & = K_{X_s}.X_s \\ &=(K_{S}+X_s).X_s \\ & = K_S.X_s + X_s^2\\ &=  rc_1(\mathcal{O}(-D)) +r^2 X^2 \\&= -rn+ r^2(2g-2 +n)
 \end{align*}
 and hence \begin{equation}\label{dimspec}g(X_s)=r^2(g-1)+\frac{rn(r-1)}{2}+1.\end{equation}

The eigenvalues of $\Phi|_x$ for $x\in X$ control the ramification of $X_s$ over $x$, e.g. if all eigenvalues are $0$, then $X_s$ is completely ramified over $x$, if all are different then $X_s$ is unramified over $x$.  We denote the covering by $\rho:X_s\rightarrow X$, with ramification divisor $R$ on $X_s$.

\subsubsection{Generic fibers}
Now, if $X_s$ is smooth, pull back $E$ to $X_s$ by $\rho$.  We canonically get a line bundle $L$ on $X_s$, such that $L(-R)$ sits inside this pull back (see e.g.~\cite{BNR},Prop.3.6).  Away from a ramification point the fiber of $L$ is given by an eigenspace of $\rho^*E$, exactly corresponding to the eigenvalue of $\Phi$ given by that point of $X_s$, in fact one has that
\begin{equation}\label{bnr} 0\rightarrow L(-R)\rightarrow\rho^*E\stackrel{\rho^*\Phi-\lambda\id}{\longrightarrow} \rho^*(E\otimes K(D) )\rightarrow L\otimes  \rho^*K(D)\rightarrow 0\end{equation} is exact, see \cite{BNR},Remark 3.7.  Furthermore we have $\rho_*L=E$, and the Higgs field $\Phi$ is also easily recovered: multiplication by the canonical section $\lambda$ of $\rho^*K(D)$ descends to a morphism $$\Phi: \rho_*L\rightarrow \rho_*L\otimes K(D) = \rho_* \left( L\otimes \rho^*K(D) \right).$$

 In order to obtain the degree of $L$, apply Grothendieck-Riemann-Roch to the morphism $\rho:X_s\rightarrow X$, and then integrate both sides.  This gives $$\deg(E)+r\frac{1}{2}\deg(T_X) = \deg (L)+\frac{1}{2}\deg (T_{X_s})$$
and hence \begin{align*}\deg(L) &= \deg(E)+r\frac{1}{2}\deg(T_X) -\frac{1}{2}\deg(T_{X_S}) \\
&= d+  r (1-g) + r^2(g-1) + \frac{rn(r-1)}{2} \\ & = d+ r(1-r)(1-g-\frac{n}{2}).\end{align*}

Moreover, the smoothness of the spectral curve guarantees that only finitely many parabolic structures are compatible with the Higgs field, even if the Higgs field has repeated eigenvalues at the points of $D$.  
Indeed, we have
\begin{proposition}\label{unique}If the spectral curve $X_s$ is smooth, then $h_{\alpha}^{-1}(s)$ consists of disjoint copies of the Jacobian of $X_s$, one for every partition of the eigenvalues along the multiplicities of the flagtype. 
\end{proposition} 
With `partition along the multiplicities' we mean here a partition of the set of eigenvalues into subsets, the sizes of which are the various multiplicities.
\begin{proof} We want to show that with each line bundle $L$ on $X_s$ of degree $d+r(1-r)(1-g)$ we get a parabolic Higgs bundle of degree $d$.  

As stated above, the push forward $\rho_*L$ determines a vector bundle $E$ on $X$~\cite{BNR}, and   multiplication by $\lambda$ on $L$ 
descends to a Higgs field.
What remains is to construct a flag of the desired type at the marked points.
If $p$ is in $D$, denote the eigenvalues of $\Phi|_p$ by $\sigma_1, \dots, \sigma_r$.
Choose a partition of the $\sigma_i$ according to the multiplicities $m_{i}(p)$ and relabel the $\sigma_i$ such that the partition is given by: 
$$\{\sigma_{1},\dots,\sigma_{ m_{1}(p)}\},\{\sigma_{m_{1}(p)+1},\dots,\sigma_{m_{1}(p)+m_{2}(p)} \},\dots,
\{\sigma_{r-m_{r(p)}} ,\dots, \sigma_r \}.$$
Now, choose a Zariski open set $W$ around $p$ such that $L$ is trivial over its inverse image $\rho^{-1}(W)$, and $K(D)$ is trivial over $W$.   Observe that we can always do this: take any rational section of $L$, then by the independence of valuation theorem (see e.g. \cite{deuring}, page 19) we can choose a rational function such that the product of the two has no poles or zeroes on $\rho^{-1}(p)$.  This new rational section trivializes $L$ on the complement of its divisor, which then easily gives the desired $W$.
Restricting to $W$, $\rho_*L$ is given as an $\mathcal{O}(W)$-module by $$\mathcal{O}(W)[x]/(x^r+s_1x^{r-1}+\dots+s_r),$$ with the $s_i\in \mathcal{O}(W)$.  Because of the choices made, we also know that at $p$, the $s_j$ are given by the elementary symmetric polynomials in the $\sigma_{i}$.  Therefore the fiber of $\rho_*(L)$ over $p$ is exactly given by $$\mathbb{C}[x]/( (x-\sigma_{1})\dots(x-\sigma_{r})).$$  $\Phi$ is of course given by multiplication by $x$ here, and hence this determines a basis: if we put $e_j=\prod_{l=1}^{j-1}(x-\sigma_{l})$, with $e_1=1$, then with respect to this basis $\Phi$ looks like the lower-triangular matrix 
\begin{equation}\label{matrix}\left(\begin{array}{ccccc} \sigma_{1} & & & & \\ 1 & \sigma_{2} & & & \\  & 1 & \ddots & & \\  & & \ddots & \ddots & \\& & & 1 & \sigma_{r}\end{array}\right).\end{equation}
Hence by using $e_1, \dots, e_r$ as an adapted basis, i.e. $$E_{p,1}=<e_1,\dots,e_r>,\ \cdots,\  E_{p,r(p)}=<e_{r-m_{r(p)}(p)},\dots, e_r>,$$ we get the parabolic structure, which is independent of our trivializations of $L$ and $K(D)$ and the choice of the basis $e_i$.  From the matrix form (\ref{matrix}) one can immediately see that the eigenspace for every eigenvalue is one-dimensional, even if the eigenvalue has multiplicity, this is what gives the uniqueness of the filtration.

It remains to show that this parabolic Higgs bundle is stable: for this observe that a $\Phi$-preserved subbundle of $E$ would necessarily correspond to a subsheaf of $L$ (see again~\cite{BNR}, page 174), and as $X_s$ is assumed to be smooth, this has to be a locally free sheaf itself as well, necessarily of lower degree.  Notice that this shows that smoothness of the spectral curve implies that there are no subbundles preserved by the Higgs field, and hence the slope-stability condition need not even be applied.
\end{proof}
In section~\ref{leaves} we shall see that each of these Jacobians is actually contained in a different symplectic leaf for the Poisson structure on $\phig$.  In section \ref{G-S} we shall further study the forgetfull morphisms,  mentioned above in (\ref{laatste}), from moduli spaces $\mathcal{P}_{\tilde{\alpha}}$ of finer flag type (e.g. full flag type) to moduli spaces $\phig$ of coarser flag type (e.g. $\mathcal{P}_0$ for the minimal flag types, with all flags being $E_{p,1}\supset E_{p,2}=\{0\}$).  Note that from theorem \ref{unique} we can already conclude that such a morphism will be finite over the locus in $\mathcal{H}$ corresponding to smooth spectral curves.  

The moduli-space $\mathcal{M}_{\text{Higgs}}$ of non-parabolic Higgs bundles (with Higgs fields $\Phi:E\rightarrow E\otimes K$) is a subvariety of $\mathcal{P}_0$; in fact it is a symplectic leaf for the Poisson structure we shall exhibit on $\mathcal{P}_0$.  All in all this gives us the following diagram (with $h$ the Hitchin map for $\mathcal{M}_{\text{Higgs}}$):
$$\xymatrix{& \mathcal{P}_{\tilde{\alpha}} \ar[d]\ar^{h_{\tilde{\alpha}}}[ddr] & \\ \mathcal{M}_{\text{Higgs}}\ar@{^{(}->}[r] \ar_h[dr] & \mathcal{P}_0 \ar|{h_{\alpha}}[dr] & \\ & \bigoplus_i H^0(X, K^i)\ar@{^{(}->}[r] &\mathcal{H}.}$$

\section{Poisson structure}\label{poissonstruc}

\subsection{Background material}
We will begin by briefly reviewing the background on Poisson geometry and Lie algebroids that we need. 
Remark that in the literature Lie algebroids and groupoids are usually used in a differential geometric setting.  For our purposes, where we will use these notions on the moduli space of parabolic Higgs bundles, all structures (bundles, spaces, actions) are algebraic however.
\\

There are two main differences in the holomorphic or algebraic setting vs. the smooth settings: algebroids have to be defined using the sheaf of sections of the underlying vector bundle rather than just the global sections, and more significantly, principal bundles do not always have connections, or equivalently the corresponding Atiyah sequence (\ref{atisequence2}) does not always split.  As we only use Lie groupoids and algebroids in a non-singular (algebraic) setting, we have however kept the differential geometric notions of submersion etc.  For more background on Lie groupoids and algebroids see~\cite{cannasweinstein,mackenzie}, which we use without reference in this section. 

\subsubsection{Poisson geometry}
There are many ways of packaging a Poisson structure, so just to fix conventions we shall state the one most convenient for our purposes:
\begin{definition}A \emph{Poisson structure} on a complex manifold $M$ is given by a bundle morphism $\sharp:T^*M\rightarrow TM$ that is anti-symmetric, i.e. $\sharp^*=-\sharp$, such that the Schouten-Nijenhuis bracket $[\tilde{\sharp},\tilde{\sharp}]$ of the corresponding bivector $\tilde{\sharp}\in \bigwedge^2 TM$ is zero.
\end{definition}
Now, for any manifold $N$, look at the total space of the cotangent bundle $\pi:T^*N\rightarrow N$.  The tangent and co-tangent bundles to $T^*N$ both fit in short exact sequences,
\begin{equation}\label{shorteen}0\rightarrow \pi^*(T^*N) \rightarrow T^*(T^*N)\rightarrow \pi^* TN \rightarrow 0\end{equation}
and 
\begin{equation}\label{shorttwee}0\rightarrow \pi^*(T^*N) \rightarrow T(T^*N)\rightarrow \pi^* TN \rightarrow 0.\end{equation}
The canonical Poisson structure (which is of course even a symplectic structure) is defined as the unique anti-symmetric bundle morphism $\sharp_{T^*N}:T^*(T^*N)\rightarrow T(T^*N)$ such that in the diagram
\begin{equation}\label{antisym}\xymatrix{ 0 \ar[r] & \pi^*(T^*N)\ar[d]_{Id}\ar[r] & T^*(T^*N) \ar[d]^{\sharp_{T^*N}} \ar[r] & \pi^* TN \ar[d]_{-Id}\ar[r] &0 \\
 0 \ar[r] & \pi^*(T^*N)\ar[r] & T(T^*N) \ar[r] & \pi^* TN \ar[r] &0}\end{equation} all squares commute.   
\begin{definition} A morphism $f:M_1\rightarrow M_2$ between Poisson spaces $(M_1, \sharp_1)$ and $(M_2,\sharp_2)$ is Poisson if the square \begin{equation}\label{antisym2}\xymatrix{ T^*M_1\ar[r]^{\sharp_1} & TM_1\ar[d]^{df}\\ T^*M_2 \ar[u]^{(df)^*} \ar[r]_{\sharp_2} & TM_2 }\end{equation} commutes. 
\end{definition}
In particular, if a Poisson manifold $M$ is equipped with an action by a group $G$ that preserves the Poisson structure, then if the quotient exists the quotient map $M\rightarrow M/G$ is a Poisson morphism.
\subsubsection{Algebroids and Poisson structures}
Two places where Poisson structures naturally occur are on the quotients of a symplectic manifold by a Hamiltonian group action - where the symplectic leaves are given by the various symplectic reductions.
Another place is on the total space of the dual of a Lie algebroid (which includes as a special example the dual of a Lie algebra).  We will mainly be interested in the special case of an Atiyah algebroid, which is an example of both of these situations. 

\begin{definition} A \emph{Lie algebroid} over a (complex) variety $M$ is a vector bundle $E\rightarrow M$ such that the sheaf of sections of $E$ is a sheaf of Lie algebras for a bracket $$[.,.]:\mathcal{O}(E)(U)\times \mathcal{O}(E)(U)\rightarrow \mathcal{O}(E)(U),$$ together with a bundle map, the \emph{anchor},  $a:E\rightarrow TM$ which preserves the Lie brackets on sections.  Moreover the following Leibniz rule has to hold, for $f\in \mathcal{O}(U), X, Y\in \mathcal{O}(E)(U):$
$$[X,fY]= f[X,Y] +  (a(X)f)Y.$$ 
\end{definition}

Two natural classes of examples of Lie algebroids are given by tangent bundles $TM$ of a manifold $M$ (where the map $a$ is the identity), and Lie algebras $\mathfrak{g}$, regarded as a vector bundle over a point.  One can think of transitive Lie algebroids (i.e. algebroids with surjective anchor maps) as interpolations between these two.

The relevance of Lie algebroids for us is through the following theorem, which in this form is due to Courant~\cite{courant} (see also~\cite{mackenzie,cannasweinstein} or~\cite{bottac}):
\begin{theorem}\label{algtopois} The total space of the dual vector bundle $E^*$ of a Lie algebroid $E$ has a natural Poisson structure.
\end{theorem}

For the two previous examples mentioned above the Poisson structures are given by the canonical symplectic structure on the total space of $T^*M$, and the Kirillov-Kostant-Souriau Poisson structure on $\mathfrak{g}^*$.  In the former there is one single symplectic leaf, in the latter case the symplectic leaves are given by the co-adjoint orbits of the Lie group $G$.
\\

Let a group $G$ act freely and properly on a manifold $P$, in other words $P\overset{\pi}{\rightarrow} P/G$ is a $G$-principal bundle.  Of course then $G$ also acts freely and properly, in a Hamiltonian fashion, on the symplectic manifold $T^*P$, and therefore the quotient $T^*P/G$ is a Poisson manifold. Another way to realize the Poisson structure on $T^*P/G$ is as the dual of a particular type of Lie algebroid, 
the so-called \emph{Atiyah algebroid}, 
as follows.  $G$ acts freely on $TP$, and $TP/G$ is a Lie algebroid over $P/G$.  One sees this most easily by interpreting the sections of $TP/G$ as $G$-invariant vector fields on $P$, and the sections of $T(P/G)$ as $G$-invariant  sections of the bundle on $P$ that is the quotient of $TP$ by the tangent spaces to the orbits: \begin{equation}\label{atisequence}0\rightarrow T_{\text{orbits}}P \rightarrow TP \rightarrow \pi^*T(P/G)\rightarrow 0.\end{equation}  The anchor map $a:TP/G\rightarrow T(P/G)$ is just given by projecting an invariant vector field to the part `orthogonal' to the orbits.  It is clear that this satisfies the required property, since functions on $P/G$ correspond to $G$-invariant functions on $P$, and any tangent field along the orbits annihilates an invariant function. 
\\

All in all the above shows that $TP/G$ is the extension of $T(P/G)$ by the adjoint bundle $\text{Ad}(P)=P\times_{\text{Ad}} \mathfrak{g}$.  The latter is, as a bundle of Lie algebras, of course a Lie algebroid with trivial anchor map, and the corresponding \emph{Atiyah sequence} \begin{equation}\label{atisequence2}0\rightarrow \text{Ad}(P)\rightarrow TP/G \rightarrow T(P/G)\rightarrow 0\end{equation} preserves all Lie brackets on local sections.
\\

Though mainly used in a differential geometric setting, the Atiyah algebroid and the corresponding short exact sequence were originally introduced in~\cite{aticonn} in the context of the study of the existence of holomorphic connections in complex fiber bundles. 

\subsubsection{Symplectic leaves for the dual of an algebroid}  
Since we are concerned with Atiyah algebroids we can study the symplectic leaves on its dual fairly directly, using the general fact that the symplectic leaves of the Poisson reduction of a symplectic manifold by a free Hamiltonian group action correspond to the various symplectic reductions.  We shall do this in section~\ref{leaves}.
For the sake of completeness we do briefly indicate here however that these symplectic leaves can be seen as co-adjoint orbits of a Lie groupoid.  For us this is mainly of philosophical relevance, leading to the interpretation of $\phig$ in the case of full flags as a Grothendieck-Springer resolution.
\\
 
Just as Lie algebras $\mathfrak{g}$ are given as tangent spaces to Lie groups $G$, Lie algebroids can come from a Lie groupoid - though not every algebroid integrates to a Lie groupoid, see~\cite{intlie}.  

\begin{definition}
A \emph{Lie groupoid} $G\rightrightarrows M$ over a manifold $M$ is a space\footnote{In the differential geometric setting it is in general not required here that $G$ is Hausdorff, but it is assumed that $M$ is.} $G$ together with two submersions $\alpha, \beta:G\rightarrow M$, as well as an associative product $(g_1,g_2)\mapsto g_1g_2 $ defined on \emph{composable pairs}, i.e. $(g_1, g_2)$ s.t. $\beta(g_2)=\alpha(g_1)$, such that $\alpha ( g_1 g_2) =\alpha (g_2)$ and $\beta  ( g_1 g_2)=\beta (g_1)$.  Furthermore an identity section $\epsilon: M\rightarrow G$ has to be given, such that the following hold for all $g\in G$: $$\epsilon(\beta (g))g=g\text{ and }g\epsilon(\alpha (g))= g$$ as well as an inversion $\iota:G\rightarrow G$, with $$\iota(g)g=\epsilon(\alpha(g)) \text{ and } g\iota(g)=\epsilon(\beta( g)).$$
\end{definition}
The maps $\alpha$ and $\beta$ are often referred to as respectively the \emph{source} and \emph{target} maps of the groupoid, and one thinks of $G$ as consisting of arrows $g$ from $\alpha(g)$ to $\beta(g)$, and compositions, inverses and identities can be understood as such.
\\

With every Lie groupoid one can naturally associate a Lie algebroid structure on the normal bundle to $M\cong \epsilon (M) \subset G$.
We refer to~\cite{cannasweinstein,mackenzie} for  
further background material. 
\\

Now, let $G\rightrightarrows M$ be a Lie groupoid, with associated Lie algebroid $E\rightarrow M$.  Then there exists a \emph{cotangent groupoid} $T^*G \rightrightarrows E^*$, 
which is both a vector bundle over $G$ and a Lie groupoid over $E^*$.
Clearly any groupoid $H\rightrightarrows N$ acts on its base $N$.  
Of particular relevance for us is the following (see e.g.~\cite{mackenzie}, Proposition 11.5.4 and Theorem 11.5.18):
\begin{theorem} The symplectic leaves for the Poisson structure on the total space of the dual of a Lie algebroid $E\rightarrow M$ associated with a Lie groupoid $G\rightrightarrows M$ are the (connected components of) the orbits for the action of $T^*G \rightrightarrows E^*$ on $E^*$.
\end{theorem}
The orbits of $T^*G\rightrightarrows E^*$ are often referred to as the co-adjoint orbits of the original groupoid $G\rightrightarrows M$.
\\

This theorem establishes a cotangent groupoid as a particular case of the general notion of a \emph{symplectic groupoid} (see~\cite{groupoidsymp}).  The base of a symplectic groupoid is always Poisson, and its symplectic leaves are given by the orbits of the symplectic groupoid~\cite{groupoidsymp}. 
\\

\subsection{Poisson structure on \phig}
\subsubsection{The $\sharp$ map} 
The tangent space to $\phig$ at a stable parabolic Higgs bundle $(E,\Phi)$ is given by the $\mathbb{H}^1$ hypercohomology of the two-term 
 complex \begin{equation}\label{een}\parend (E) \overset{[.\, ,\Phi]}{\longrightarrow} \parend(E)\otimes K(D).\end{equation}
 
Let us now write down the Poisson bracket.  The dual of the complex (\ref{een}), tensored with $K$, is given by  
\begin{equation}\label{anderhalf}\sparend (E) \overset{-[.\, ,\Phi]}{\longrightarrow} \sparend(E)\otimes K(D).\end{equation} We can now inject (\ref{anderhalf}) into (\ref{een}), as follows:
\begin{equation}\label{tekens}
\xymatrix{ 
\sparend{E}\ar_{-[.\,,\Phi]}[d] \ar@{^{(}->}^{\id}[r] 
&\parend (E)\ar^{[.\,,\Phi]}[d] \\ 
\sparend(E)\otimes K(D)\ar@{^{(}->}^{-\id \otimes \id_{K(D)}}[r]
& \parend(E)\otimes K(D).
}\end{equation}
  Using Serre duality for hypercohomology we therefore get a map 
\begin{multline}\label{sharpie}\sharp_{\phig}:T^*_{[E,\Phi]}\phig \cong \mathbb{H}^1\left( \sparend (E) \rightarrow \sparend(E)\otimes K(D)\right)\\ \rightarrow \mathbb{H}^1\left( \parend (E) \rightarrow \parend(E)\otimes K(D)\right)\cong T_{[E,\Phi]}\phig.\end{multline} 
Because of the choice of signs in (\ref{tekens}) $\sharp_{\phig}$ is antisymmetric.  We now want to show that this determines a Poisson structure on $\phig$\footnote{We would like to point out the similarity between the bivector (\ref{sharpie}) and the bivector obtained by Bottacin in his study of Poisson structures on moduli spaces of parabolic vector bundles on algebraic surfaces, see \cite{bottac2}, equation 4.1.  As Higgs bundles originally arose through a dimensional reduction of a structure in higher dimensions this is not surprising.}. 
One could try to do this directly by calculating the Schouten-Nijenhuis bracket, but it would be rather hard and not so instructive, therefore we follow a different route below.

\subsubsection{Poisson structure via Lie algebroids}
Now let $\phig^0$ be the open subvariety of $\phig$ consisting of those parabolic Higgs bundles $(E,\Phi)$ whose underlying parabolic bundle is stable.  As mentioned before, this is a vector bundle over  $\pvect$, the moduli space of parabolic vector bundles, with fiber $H^0(X, \parend(E)\otimes K(D) )$.

We shall use the following two projections:
\begin{equation}\label{projections}
\xymatrix{& \pvect \times X\ar[dl]_{\eta} \ar[dr]^{\nu} & \\ \pvect & & X} \end{equation}

Now, on $\pvect\times X$ we have a parabolic universal bundle~\cite[Theorem 3.2]{bodenyoko2}, which we denote by $\mathcal{E}$.  This leads to a short exact sequence of sheaves on $\pvect\times X$:
$$0\rightarrow \sparend (\mathcal{E})\rightarrow \parend (\mathcal{E}) \rightarrow \prod_{p\in D} \mathfrak{l}_p\otimes\mathcal{O}_{{\nu^{-1}}(p)}\rightarrow 0.$$
Applying ${\eta_*}$ to this sequence gives the exact sequence
\begin{multline}\label{long}
0\rightarrow 
{\eta_*}\sparend (\mathcal{E})
\rightarrow {\eta_*} \parend (\mathcal{E})
 \rightarrow {\eta_*}\left(\prod \mathfrak{l}_p\otimes\mathcal{O}_{{\nu^{-1}}(p)}\right) \rightarrow
  \\ 
  R^1{\eta_*}\sparend (\mathcal{E}
 )
 \rightarrow R^1{\eta_*} \parend (\mathcal{E}) 
 \rightarrow R^1{\eta_*}\left(\prod \mathfrak{l}_p\otimes\mathcal{O}_{{\nu^{-1}}(p)}\right) 
 \rightarrow 0.
\end{multline}
As the support of $\prod \mathfrak{l}_p\otimes \mathcal{O}_{\nu^{-1}(p)}$ has relative dimension zero with respect to $\eta$, the last term of this sequence is easily seen to be zero by relative dimension vanishing (see e.g. \cite{hartshorne}, III.11.2). The first term is zero since every stable bundle is simple (lemma~\ref{stablesimple}), and for the same reason the second term\footnote{Observe that this term would also be zero in the case of a semi-simple structure group, so for parabolic principal Higgs bundles the relevant principal bundle over the moduli space $\pvect$ would just have the product of the Levi groups as structure group, without quotienting by global endomorphisms.}
 is an invertible sheaf. 
 Moreover, $R^1{\eta_*} \parend (\mathcal{E}) $ is the tangent sheaf to $\pvect$, hence we shall denote it as $\mathcal{T}_{\pvect}$.
\begin{lemma}All of the sheaves occurring in the sequence (\ref{long}) are locally free.
\end{lemma}
\begin{proof} It suffices to notice that the corresponding cohomology groups (e.g. $H^1(\sparend(E) )$) have constant rank as $E$ varies in $\pvect$, and apply  Grauert's theorem (\cite{hartshorne}, III.11.2).
\end{proof}

Denote now ${\eta_*}\left(\prod \mathfrak{l}_p\otimes\mathcal{O}_{{\nu^{-1}}(p)}\right) / {\eta_*} \parend (\mathcal{E})$ by $\text{Ad}$.  This clearly is a bundle of Lie algebras.  Our claim is that the short exact sequence \begin{equation}\label{atisequence3}
0\rightarrow \text{Ad} \rightarrow R^1{\eta_*} \sparend(\mathcal{E}) \rightarrow \mathcal{T}_{\pvect}\rightarrow 0\end{equation} is an Atiyah sequence, as in (\ref{atisequence2}).

\begin{lemma}The dual of the vector bundle $\phig^0 \rightarrow \pvect$ is $R^1{\eta_*}\sparend(\mathcal{E})$.
\end{lemma}
\begin{proof}By relative Serre duality for the morphism $\eta$ in (\ref{projections}), the dual of $R^1{\eta_*}\sparend(\mathcal{E})$ is given by ${\eta_*}\parend(\mathcal{E}\otimes \nu^*K(D))$, which is a locally free sheaf as well, with fibers of the associated vector bundle over a point $[E]\in\pvect$ given by $H^0(\parend(E\otimes K(D) ))$.  Clearly the obvious bundle morphism ${\eta_*}(\sparend(\mathcal{E}\otimes \nu^* K(D) )\rightarrow \phig^0$ is an isomorphism.
\end{proof}

We now assume the existence of a principal bundle $\pi:\pfram\rightarrow \pvect$ with structure group \begin{equation}\label{struct}\mathcal{L}= \left( \prod_{p\in D} L_p \right)/C, \end{equation} where $C$ is the diagonal subgroup of the product of the centers of the $P_p$, such that
the total space of $\pfram$ can be interpreted as a moduli space of framed $\alpha$-stable parabolic bundles. A framed parabolic bundle here is a parabolic vector bundle together with a \emph{framing} of the associated graded space at the marked points, i.e. the choice of an isomorphism $$\Gr (p)=\bigoplus_{i=1}^{r(p)} F_{p,i}/F_{p,i+1}\overset{\cong}{\longrightarrow}\bigoplus \mathbb{C}^{m_i(p)} $$ at each marked point $p\in D$.  In appendix~\ref{append} a construction for $\pfram$ is given in the case of full flags.

\begin{theorem}\label{tired} The sequence (\ref{atisequence3}) is the Atiyah sequence for the $\mathcal{L}$-principal bundle $\pfram\rightarrow \pvect$.
\end{theorem}
\begin{proof}  
Since the sheaf of sections of the adjoint bundle of a principal bundle is the direct image of the relative tangent sheaf of the associated projection, it suffices to show that the inverse image under $\pi$ of the short exact sequence $(\ref{atisequence})$ is the sequence on $\pfram$ determining the relative tangent sheaf: $$0\rightarrow \mathcal{T}_{\pi}\rightarrow \mathcal{T}_{\pfram}\rightarrow \pi^*\mathcal{T}_{\pvect}\rightarrow 0 .$$  To fix notation, let us look at the commutative diagram \begin{equation}\label{commu}\xymatrix{ 
\pfram\times X 
\ar[r]^{\widetilde{\pi}} \ar[d]_{\widetilde{\eta}} 
& \pvect\times X \ar[d]^{\eta} \\ 
\pfram
\ar[r]^{\pi} 
&\pvect}
\end{equation}
where $\pfram$ is a moduli space of framed stable parabolic bundles.  The space $\pfram\times X$ comes equipped with a sheaf of framed endomorphisms (endomorphisms that preserve the framing) of a universal bundle $\fparend$, such that the tangent sheaf to $\pfram$ is given by $R^1{\widetilde{\eta}}_* \fparend$.  As it is also easy to see that $\fparend$ is equal to $\widetilde{\pi}^* \sparend(\mathcal{E})$, the commutativity of the diagram (\ref{commu}) - and flatness of $\pi$, guaranteed by invoking e.g.~\cite{grothendieck}, proposition 6.1.5  - give that indeed $$R^1{\widetilde{\eta}}_* \fparend\cong \pi^*R^1{\eta_*}\sparend(\mathcal{E}),$$ see e.g. \cite{hartshorne}, Proposition III.9.3.
\end{proof}

\subsubsection{Extension of the bracket}

\begin{theorem}\label{endinsight}
The bivector on $\phig$ determined by (\ref{sharpie}) extends the Poisson structure on $\phig^0$ given by the Lie algebroid on $(\phig^0)^*$. 
\end{theorem} 
\begin{proof}
For an  $\alpha$-stable framed parabolic bundle $(E,\cong)$ the tangent and cotangent spaces to $\pfram$ are given respectively by $$T_{[E,\cong]}\pfram =H^1(\sparend(E))$$ and  $$T^*_{[E,\cong]} \pfram = H^0(\parend(E)\otimes K(D)).$$
The tangent and the cotangent spaces to $T^*\pfram$ at a point $[E,\Phi,\cong]$ are given by the first hypercohomology groups $$T_{[E,\Phi,\cong]} T^*\pfram\cong \mathbb{H}^1\left({\scriptstyle \sparend(E)\overset{[.\, ,\Phi]}{\longrightarrow}\parend{E}\otimes K(D)} \right)$$ 
and 
$$T^*_{[E,\Phi,\cong]} T^*\pfram \cong\mathbb{H}^1\left({\scriptstyle \sparend(E)\overset{-[.\, ,\Phi]}{\longrightarrow}\parend{E}\otimes K(D)} \right).$$
Furthermore, if we look at the short exact sequences of complexes
$$\xymatrix{ 0\ar[r] & 0 \ar[r]\ar[d] & \sparend(E)\ar[r]\ar_{\pm[,\, ,\Phi]}[d] & \sparend(E) \ar[r]\ar[d] & 0\\ 0 \ar[r]& \parend(E)\otimes K(D)\ar[r] & \parend(E)\otimes K(D)\ar[r]
  & 0 \ar[r] & 0
 },$$
take their long exact sequence in hypercohomology, and notice that\\ $H^0(\sparend(E))$ is zero since $E$ is stable as a parabolic vector bundle, we obtain the short exact sequences
\begin{multline*}
 0\rightarrow H^0(\parend(E)\otimes K(D))\rightarrow \mathbb{H}^1\left(
 {{\tiny{\begin{array}{c} \sparend(E) \\ \pm [.\, ,\Phi] \downarrow \\ \parend (E) \otimes K(D) \end{array}}}} 
 \right)\\ \rightarrow H^1(\sparend(E))
 \rightarrow 0,\end{multline*} corresponding to (\ref{shorteen}) and (\ref{shorttwee}).
Now, using the characterization given in (\ref{antisym}) for the canonical Poisson structure on $T^*\pfram$, one can see that this is induced by the morphism of complexes 
\begin{equation*}
\xymatrix{
 \sparend{E}\ar_{-[.\, ,\Phi]}[d] \ar@{^{(}->}^{\id}[r] 
 &\sparend (E) \ar_{[.\, ,\Phi]}[d] 
 \\ \parend(E)\otimes K(D)\ar@{^{(}->}^{-\id \otimes \id_{K(D)}}[r] 
 & \parend(E)\otimes K(D).
}.\end{equation*}
Since the map $T^*\pfram \rightarrow \phig^0$ is a Poisson morphism, and using the characterization (\ref{antisym2}) and the definition (\ref{sharpie}) of $\sharp_{\phig}$ it suffices to notice that by these choices indeed the square 
$$\xymatrix{
\mathbb{H}^1\left({{\tiny \begin{array}{c} \sparend(E) \\ {- [.\, ,\Phi]} \downarrow \\ \parend (E) \otimes K(D) \end{array}}}\right) \ar[r] &  
\mathbb{H}^1\left( {\tiny \begin{array}{c} \sparend(E) \\ { [.\, ,\Phi]} \downarrow \\ \parend (E) \otimes K(D) \end{array}} \right) \ar[d] \\
\mathbb{H}^1\left({\tiny \begin{array}{c} \sparend(E) \\  {- [.\, ,\Phi]} \downarrow \\ \sparend (E) \otimes K(D) \end{array}} \right) \ar[r]\ar[u] & 
\mathbb{H}^1\left({\tiny \begin{array}{c} \parend(E) \\   { [. \, ,\Phi]} \downarrow \\ \parend (E) \otimes K(D) \end{array}} \right)
}$$

commutes.
\end{proof}
 
 As $\phig^0$ is open and dense in $\phig$ this establishes the Poisson structure on all of $\phig$.

\subsubsection{Symplectic leaves}\label{leaves}
Given a principal $G$-bundle $P\rightarrow M$, the symplectic leaves of the dual of an Atiyah algebroid $TP/G$ are simply the symplectic reductions of the cotangent bundle $T^*P$.  The lift of the action of a group $G$ on a manifold $N$ to the total space of the cotangent bundle $T^*N$ is of course always Hamiltonian, with a canonical moment map given by $$g\in \mathfrak{g},\ \chi\in T^*_x N\  :   \ \ \ \mu(\chi)(g)=\chi(\xi_g(x)), $$ 
where $\xi_g$ is the Hamiltonian vector field corresponding to $g\in\mathfrak{g}$.  In the particular case where the action of $G$ is free, i.e. when the manifold is a principal $G$-bundle, the moment map can also be understood by dualizing the sequence (\ref{atisequence}):
$$0\rightarrow\pi^*T^*(P/G)\rightarrow T^*P\overset{\mu}{\rightarrow} T^*_{\text{orbits}}P\rightarrow 0,$$
observing that $T_{\text{orbits}}P\cong P\times \mathfrak{g}$ and $T^*_{\text{orbits}}P\cong P\times\mathfrak{g}^*$.  Recalling (\ref{long}) and the proof of Theorem~\ref{tired} this tells us immediately that on $T^*\pfram$ the moment map $\mu$ is given by the \emph{parabolic residue}, i.e. the map that, when restricted to a fiber over a framed parabolic bundle $E$, gives
$$H^0(\parend(E)\otimes K(D))\rightarrow \ker \left( \bigoplus_{p\in D} \mathfrak{l}^*_p\rightarrow H^0(\parend(E))^*\right)$$ which is the dual of the boundary map associated with the short exact sequence $$0\rightarrow\sparend(E)\rightarrow\parend(E)\rightarrow \bigoplus_{p\in D}\mathfrak{l}_p\otimes \mathcal{O}_p\rightarrow 0.$$

The reduction, and hence symplectic leaves, are simply $\mu^{-1}(O)/\mathcal{L}$, where $O$ is a co-adjoint orbit in $\text{Lie}^*(\mathcal{L})$.  Notice that in the generic case, when the eigenvalues of the Higgs field at the marked points are all distinct, there is a unique co-adjoint orbit with these eigenvalues.  When eigenvalues are repeated on a particular $E_{p,i}/E_{p,i+1}$ there will be several co-adjoint orbits.
\\
 
 One can check that this agrees with the rank of the Poisson structure given by $\sharp_{\phig}$ at a parabolic Higgs bundle $(E,\Phi)$.   Indeed, if one recalls the definition (\ref{sharpie}) then the short exact sequence of complexes 
 $$\xymatrix{ 0 \ar[r]& \sparend(E)\ar^{\id}[r] \ar_{-[.\, ,\Phi]}[d] & \parend(E)\ar_{ [. \, ,\Phi]}[d]\ar[r] & \bigoplus \mathfrak{l}_p\otimes \mathcal{O}_p\ar[r] \ar^{ [.\,  ,\Phi]|_{\Gr (p)}}[d] & 0 \\ 
 0\ar[r] & \sparend(E)\otimes K(D)\ar[r] & \parend(E)\otimes K(D) \ar[r] & \bigoplus \mathfrak{l}_p\otimes K(D)_p\ar[r]  & 0 
 }$$
 gives rise to the long exact sequence of hypercohomology (at least when $E$ is stable):
 $$0\rightarrow\mathbb{C}\rightarrow \oplus \ker\left([.\, ,\Phi]|_{\Gr (p)}\right)\rightarrow \mathbb{H}^1\left( {{\tiny{\begin{array}{c} \sparend(E) \\ \downarrow \\ \sparend (E) \otimes K(D) \end{array}}}}\right)\overset{\sharp_{\phig}}{\longrightarrow} \mathbb{H}^1\left( {{\tiny{\begin{array}{c} \parend(E) \\ \downarrow \\ \parend (E) \otimes K(D) \end{array}}}}\right) \rightarrow ... ,$$ and hence the generic, maximal rank of $\sharp_{\phig}$, occurring when all eigenvalues of $\Phi$ are different, is \begin{equation}\label{rank}\rk \phig=\dim \phig -nr+1=(2g-2)r^2 +nr(r-1) +2.\end{equation}

The eigenvalues of the Higgs field are in fact determined by the Hitchin map, and as the latter is blind to the parabolic structure, it factors through morphisms of the form (\ref{laatste}):
\begin{equation}\label{casi}\xymatrix{\mathcal{P}_{\tilde{\alpha}} \ar[d]\ar[dr]^{h_{\tilde{\alpha}}} & & \\ \phig\ar[r]_{h_{\alpha}} & \mathcal{H} \ar[r]^{e} & \mathcal{H}/\mathcal{H}_0,
} 
\end{equation}

where $\mathcal{H}_0$ is the subspace of $\mathcal{H}$ given by
$$\mathcal{H}_0=H^0(X,K)\oplus H^0(X, K^2(D))\oplus\cdots\oplus H^0(X,K^r(D^{r-1})).$$  Observe that the Hitchin map for the moduli space of strongly parabolic Higgs bundles (which is a symplectic leaf of $\phig$) takes values in $\mathcal{H}_0$.  One easily sees that, roughly speaking, $\mathcal{H}/\mathcal{H}_0$ determines the eigenvalues of the Higgs field at the marked points $p\in D$, without ordering.  In the case where $\tilde{\alpha}$ corresponds to full flags, the connected components (one for every ordering of the eigenvalues) of the fibers of $e\circ h_{\tilde{\alpha}}$ are exactly the symplectic leaves.

\subsubsection{Complete integrability of the Hitchin system}\label{integrability}
The composition $e\circ h_{\alpha}$ from (\ref{casi}) also plays a role in the Hitchin system, which we can discuss now that we have the Poisson structure at our disposal.  Recall that for a holomorphic Poisson manifold or Poisson variety of dimension $2k+l$, where the rank (or dimension of the generic leaf) of the Poisson structure is $2k$, a completely integrable system is given by $k+l$ Poisson-commuting, functionally independent functions, such that $l$ of them are \emph{Casimirs}, i.e. they Poisson-commute with any function.  Furthermore the generic fiber of the collective of these functions is required to be an abelian variety.  The connected components of the simultaneous fibers of the Casimir functions are the closures of the top-dimensional symplectic leaves.
\\

By Riemann-Roch one gets (assuming that $n\geq 1$)$$\dim(\mathcal{H})= (2g-2+n)\frac{(r+1)r}{2}+r(1-g)$$ and
$$\dim(\mathcal{H}_0)=(2g-2)\frac{r(r+1)}{2} + n \frac{r(r-1)}{2} + r(1-g)+1,$$ hence $$\dim\left(\mathcal{H}/\mathcal{H}_0\right)= nr-1.$$  Notice also from (\ref{dimspec}) that $g(X_s)=\dim\left(\mathcal{H}_0\right)$,  and from (\ref{dimyoko}) that $$2\dim\left(\mathcal{H}_0\right) + \dim\left( \mathcal{H}/\mathcal{H}_0\right) = \dim \phig.$$

In order to show directly that $\phig$ equipped with $h_{\alpha}$ is a completely integrable system, we would have to work with local information, in order to establish the vanishing of the relevant Poisson brackets.  As our map $\sharp_{\phig}$ is defined fibre-wise however, we use an alternative characterization (following \cite{markman}, section 8.1):
\begin{proposition}\label{tussen} The connected components of a generic fiber $h_{\alpha}^{-1}(s)$, corresponding to a smooth spectral curve $X_s$, are Lagrangian in symplectic leaves for the Poisson structure on $\phig$.
\end{proposition}
In order to prove this we shall need two lemmas.  The first is the following characterization of coisotropic submanifolds of symplectic leaves, which is easy to see:
\begin{lemma}\label{lemmaatje}Let $J$ be a submanifold of a symplectic leaf $L$ of a Poisson manifold $(M,\sharp)$.  Then $J$ is co-isotropic in $L$ if for any point $p\in J$ we have that $$\sharp(N^*_p J) \subset T_pJ ,$$ where $T_p J$ is the tangent space to $J$ at $p$, $N_p J$ is the normal space at $p$ of $J$ in $M$, and $N^*_p J$ is the conormal space at $p$.  
\end{lemma}
The second is a description of the tangent space to $h^{-1}_{\alpha}(s)$: 
\begin{lemma}\label{lemmaatje2}
Let $(E,\Phi)$ be a Higgs bundle in one of the components of $h_{\alpha}^{-1}(s)$, which we identify with the Jacobian $J_s$ of $X_s$.  Then we have the following short exact sequence on $X$:
$$0\rightarrow T_{(E,\Phi)}J_s \cong H^{1}(\rho_* \mathcal{O}_{X_s})\rightarrow \mathbb{H}^1\left(
 {{\tiny{\begin{array}{c} \parend(E) \\  {[.\, ,\Phi]} \downarrow \\ \parend (E) \otimes K(D) \end{array}}}} \right) \rightarrow H^0(( \rho_*K_{X_s})(D) )\rightarrow 0 $$
\end{lemma}
\begin{proof} Observe that by using Hurwitz' theorem (\cite{hartshorne}, Proposition IV.2.3) we have that $\rho^*(K_X)(R) = K_{X_s}$, and by using Hurwitz' theorem and relative Serre duality we have $\rho_*(L^{-1}(R))=E^{*}$. If we tensor the exact sequence (\ref{bnr}) on $X_s$ with $L^{-1}(R)$ and push it forward by $\rho$ we obtain the exact sequence
$$ 0\rightarrow \rho_*(\mathcal{O}_{X_s})\rightarrow \parend(E)\overset{[.\, ,\Phi]}{\longrightarrow}\parend(E)\otimes K(D)
\rightarrow \rho_*(K_{X_s})(D)\rightarrow 0$$ on $X$.  Using this we can look at the short exact sequences of complexes
$$\xymatrix{
0\ar[r] 
& \rho_*(\mathcal{O}_{X_s})\ar[d]\ar[r] 
& \parend{E} \ar_{[.\, , \Phi]}[d]\ar[r] 
& \im([.\, , \Phi])\ar[d]\ar[r] &0  \\ 0 \ar[r] & 0\ar[r] & \parend(E)\otimes K(D) \ar[r] & \parend(E)\otimes K(D) \ar[r] & 0 
}$$
and $$\xymatrix{0\ar[r]\ &\im([.\, , \Phi]) \ar[r]\ar[d] &\im([.\, , \Phi]) \ar[r]\ar[d] &0 \ar[r]\ar[d] &0 \\ 0\ar[r] & \im([.\, ,\Phi])\ar[r] & \parend(E)\otimes K(D) \ar[r] & \rho_*(K_{X_s})(D)\ar[r] &0
}
$$  Combining the hypercohomology long exact sequences gives the desired result.
\end{proof}
\begin{proof}[Proof of Proposition~\ref{tussen}]By the discussion in section~\ref{leaves}, it is clear that the connected components of $h_{\alpha}^{-1}(s)$ are contained in symplectic leaves.  Notice that by (\ref{dimspec}), Proposition~\ref{unique} and (\ref{rank}) we already know that generically these connected components are smooth subvarieties of half the dimension of the symplectic leaves, therefore it suffices to show that the connected components of $h^{-1}_{\alpha}(s)$ are co-isotropic.
For this we can now apply Lemma~\ref{lemmaatje} to a connected component of $h_{\alpha}^{-1}(s)$ corresponding to a smooth spectral curve $X_s$ by observing that in the commutative diagram
$$\xymatrix{
  H^0(\rho_*K_{X_s}(D))^*\cong H^1(\rho_* \mathcal{O}_{X_s} (-D)) \ar[r]\ar[d] 
 & H^1(\rho_*\mathcal{O}_{X_s})\ar[d]
  \\ \mathbb{H}^1\left({\scriptstyle \sparend(E)\overset{-[.\, ,\Phi]}{\rightarrow} \sparend(E)\otimes K(D) }\right) \ar^{\sharp_{\phig}}[r] \ar[d] 
  & \mathbb{ H}^1\left({\scriptstyle \parend(E)\overset{[.\, ,\Phi]}{\rightarrow} \parend(E)\otimes K(D) }\right) \ar[d]
  \\ H^1(\rho_*\mathcal{O}_{X_s})^*\cong H^0(\rho_*K_{X_s} )\ar[r] & H^0(\rho_*K_{X_s}(D) )    }$$ 
the columns (given by Lemma \ref{lemmaatje2}) are exact.
This ends the proof of Proposition~\ref{tussen}.\end{proof}
This finally gives us:
\begin{theorem}\label{groot} The moduli spaces $\phig$, with the Poisson structure introduced above and the Hitchin map $h_{\alpha}$, form completely integrable systems, for which the Casimirs are given by $e\circ h_{\alpha}$.
\end{theorem}

\section{Morphisms between moduli spaces, Grothendieck-Springer resolution}\label{G-S}
Given a complex semi-simple connected Lie Group $G$ with Lie algebra $\mathfrak{g}$ and Weyl group $W_G$, one can construct the so-called Grothendieck-Springer morphism.  There are various incarnations of this, for the group, the Lie algebra, etc, so we just briefly recall this here.  The Grothendieck-Springer space is defined as $$\mathcal{GS}_G=\{(g,\mathfrak{b})| g\in \mathfrak{g},\ \mathfrak{b}\in G/B,\ g\in\mathfrak{b}\},$$ where $B$ is a Borel subgroup\footnote{One can generalize this to $G/P$, that is, to parabolic subgroups other than Borels.} of $G$.
The obvious map $\mu:\mathcal{GS}_G\rightarrow\mathfrak{g}$  is widely used in geometric representation theory, see e.g.~\cite{ginzchriss}.  It is generically finite ($|W_G|:1$), and provides a resolution of singularities of the nilpotent cone $\mathfrak{n}\subset\mathfrak{g}$, which is referred to as \emph{Springer's resolution}.  After choosing an equivariant identification $\mathfrak{g}\cong\mathfrak{g^*}$ we can think of $\mathcal{GS}_G$ as the dual of an algebroid over $G/B$, and the map $\mu$ as the moment-map for the induced $G$ action.  In particular $\mathcal{GS}_G$ is a regular Poisson manifold (i.e. all the symplectic leaves have the same dimension), and $\mu$ will be a Poisson morphism.  Moreover there is the following diagram, called the \emph{Grothendieck simultaneous resolution}:
 \begin{equation}\label{simulres}\xymatrix{ \mathcal{GS}_G \ar[r]\ar[d] & \mathfrak{t}\ar[d] \\ \mathfrak{g}\ar[r] & \mathfrak{t}/{\mathbb{W}}_G
 }
 \end{equation} where $\mathfrak{t}$ is the abstract Cartan and $\mathbb{W}_G$ the abstract Weyl group of $G$.  For more details regarding this we refer to \cite{ginzchriss}, section 3.1.
\\

Our construction gives a similar picture for Atiyah algebroids rather than Lie algebras, where the role of $\mathfrak{g}$ is now played by any of the $\phig$, but in particular can be the moduli space of parabolic Higgs bundles with minimal flag-type (see also section~\ref{further} below), and the role of the Grothendieck-Springer variety by the moduli space of parabolic Higgs bundles with full flags. Indeed, we show below easily that this and similar forgetful morphisms are Poisson.
\\

Let us look at the moduli spaces for two different flag types on the same divisor of marked points, $\mathcal{P}_{\tilde{\alpha}}$  and $\phig$, where the flag type of the latter is coarser than that of the former.  We assume that the parabolic weights $\alpha$ and $\tilde{\alpha}$ are close enough that if one forgets part of the flag on an $\tilde{\alpha}$-stable parabolic Higgs bundle the result is $\alpha$-stable, so that we obtain a morphism \begin{equation}\label{morph}\mathcal{P}_{\tilde{\alpha}}\rightarrow \phig.\end{equation}
\begin{proposition}\label{oef}
The morphism (\ref{morph}) is Poisson.
\end{proposition}
\begin{proof}
Let us denote a parabolic Higgs bundle for the finer flag type as $(\widetilde{E},\widetilde{\Phi})$and its image under the forgetful morphism as $(E,\Phi)$.  Then clearly we have the natural inclusions of sheaves
$$\sparend(E)\subset \sparend(\widetilde{E})$$ and
$$\parend(\widetilde{E})\subset \parend(E).$$
Therefore we get that the diagram
$$\xymatrix{
\mathbb{H}^1\left({{\tiny \begin{array}{c} \sparend(\widetilde{E}) \\ \downarrow \\ \sparend (\widetilde{E}) \otimes K(D) \end{array}}}\right) \ar[r] &  
\mathbb{H}^1\left( {\tiny \begin{array}{c} \parend(\widetilde{E}) \\ \downarrow \\ \parend (\widetilde{E}) \otimes K(D) \end{array}} \right) \ar[d] \\
\mathbb{H}^1\left({\tiny \begin{array}{c} \sparend(E) \\ \downarrow \\ \sparend (E) \otimes K(D) \end{array}} \right) \ar[r]\ar[u] & 
\mathbb{H}^1\left({\tiny \begin{array}{c} \parend(E) \\ \downarrow \\ \parend (E) \otimes K(D) \end{array}} \right)
}$$
commutes.
\end{proof}

As said above, it is particularly interesting to look at the morphism (\ref{morph}) in the case where $\mathcal{P}_{\tilde{\alpha}}$ corresponds to full flags, as then $\mathcal{P}_{\tilde{\alpha}}$ is a regular Poisson manifold. 
We can put things together in the analogue of the Grothendieck simultaneous resolution~(\ref{simulres}):
\begin{equation}
\xymatrix{\mathcal{P}_{\tilde{\alpha}} \ar[d]\ar^o[rrrr]\ar^{h_{\tilde{\alpha}}}[dr] & & &  & \mathbb{C}^{rn}/\mathbb{C}\ar[d] \\ \phig \ar^{h_{\alpha}}[r] &\mathcal{H} \ar[r] &\mathcal{H}/\mathcal{H}_0 \ar^{\cong}[rr] &  &\left(\mathbb{C}^{rn}/(S_r)^n\right)/\mathbb{C}.
}
\end{equation}
The map $o:\mathcal{P}_{\tilde{\alpha}}\rightarrow \mathbb{C}^{rn}\rightarrow \mathbb{C}^{rn}/\mathbb{C}$ is just given by the eigenvalues of $\Phi$ at the marked points - because of the full flags they come with an ordering.

\section{Further remarks}\label{further}

\subsection{Comparison with Bottacin-Markman}

A particular case, the case of minimal flags, of the above has already been discussed in the literature, in independent work by Bottacin~\cite{bottac} and Markman~\cite{markman}, though not framed in terms of parabolic (Higgs) bundles or algebroids.  Reviews of this work also appeared in \cite{donagimark, donagirev}.  Bottacin and Markman study stable pairs or twisted Higgs bundles, i.e. a vector bundle $E$ over a curve $X$ together with a morphism $\Phi: E\rightarrow E\otimes F$, for some fixed line bundle $F$.  Pairs of this kind (working over a field of positive characteristic), their moduli stack and the Hitchin fibration for them also played a crucial role in the recent work of Laumon-Ng\^o~\cite{laumngo}.  A moduli space for these was constructed by Nitsure in \cite{nitsure}, and in~\cite{bottac} and~\cite{markman} it is shown that, if $\deg(F)>\deg(K)$ (or $F=K$) and once one chooses an effective divisor $D$ in $FK^{-1}$, this space has a canonical Poisson structure.  
\\

Once this choice is made, and if $D$ is moreover reduced, such a stable pair $(E, E\overset{\Phi}{\rightarrow}E\otimes K(D))$ can of course also be interpreted as a parabolic Higgs bundle for the minimal flag type $E_{p,1}\supset E_{p,2}= \{0\}$.  In the case of such minimal flags there is only a single weight at each marked point, and one sees that it does not contribute to the slope inequality (\ref{slopeineq}).  Therefore one cannot afford the luxury of the assumption of genericity of the weights, and unless the rank and degree are coprime there are properly semi-stable points, and the moduli space of vector bundles is singular.  Even over the non-singular locus, corresponding to the stable vector bundles, there does not exist a universal bundle, but there is however still a sheaf over the stable locus playing the role of sheaf of endomorphisms of a universal bundle~\cite{bottac}, Remark 1.2.3, which is as useable as our $\parend(\mathcal{E})$\footnote{Notice that in the case of parabolic bundles with non-generic weight one could use the same strategy, in fact for some of the non-generic weights an actual universal bundle does exist over the stable locus, see~\cite{bodenyoko2}, Theorem3.2.}.
\\

Both Bottacin and Markman are primarily focused on the Poisson structure, and make no mention of Lie algebroids.  Nevertheless, Bottacin uses the same philosophy of obtaining the Poisson structure through studying the dual vector bundle.  He even writes down the definition of a Lie algebroid and proves theorem~\ref{algtopois} in~\cite{bottac}, section 4.2.  He however does not identify the algebroid as an Atiyah algebroid, but rather exhibits the Lie bracket on local sections explicitly on the level of cocycles and cochains.
\\

Markman does use the principal bundle over the moduli space of bundles (using a construction of Seshadri~\cite{seshadri}), but phrases everything in terms of reduction of its cotangent bundle.  Despite this our approach is closest to Markman, and a careful reader might find several similarities in our exposition, in particular in section~\ref{integrability}, for which we were helped by \cite{markman}, section 8.1.
\\

Neither Bottacin or Markman make the restriction that we do that $D$ is a reduced divisor - i.e. they allow the Higgs field $\Phi$ to have poles of arbitrary order, when interpreted as a meromorphic bundle morphism from $E$ to $E\otimes K$, though Bottacin assumes $D$ to be reduced in some of his proofs.  Also for non-minimal flags this would be a desirable property, in particular in the light of the geometric Langlands program with wild ramification, and should not be significantly more complicated.
\\

Another obvious generalization would be to look at other semi-simple or reductive structure groups, replacing the use of spectral curves with cameral covers.  Though most of the statements we make can at least formally be translated into this setting, we have refrained from working in this generality as it seems that the dust has not settled on the notion of stability for parabolic principal bundles, cf.~\cite{bhosram,telwood,bbn1,bbn2}.  By working in the context of stacks rather than moduli schemes these problems would of course be avoided, and we intend to take up this matter in the future.

\subsection{Parabolic vs. orbifold bundles}
In the case where all weights are rational there is an alternative description of parabolic bundles in terms of \emph{orbifold} bundles,  which provided much of the original motivation (see \cite{sesori}, a history of the genesis of parabolic bundles is given in \cite{nithist}). 

Given a finite group $\Gamma$ acting on a curve $Y$, giving rise to the ramified covering $$p: Y\rightarrow X=Y/\Gamma,$$ an orbifold bundle is a $\Gamma$-equivariant bundle on $Y$.  Alternatively, in the analytic category, one can define an orbifold Riemann surface to be a (compact) Riemann surface $X$ with $n$ marked points $p_1,\ldots, p_n$ on $X$ and a positive integer $\alpha_i$ associated to each $p_i$.  An orbifold bundle is then determined by local orbifold trivialisations and transition functions, where near a marked point $p_i$ the trivializations should be of the form $D\times \mathbb{C}^r/ \sigma_i\times \tau_i$, where $D$ is a  disk in $\mathbb{C}$, $\sigma_i$ is the standard representation of $\mathbb{Z}_{\alpha_i}$, and $\tau$ is an isotropy representation $\tau:\mathbb{Z}_{\alpha_i}\rightarrow GL_r(\mathbb{C})$ $X$.  These two definitions of orbifold bundles are equivalent (under the condition that $n>2$ if $g=0$), see e.g. \cite{furutasteer}, page 42.  

An orbifold bundle in this sense corresponds to a parabolic bundle on $X$ with rational weights. The correspondence has been extended to both higher dimensions in \cite{biswasorbi} and principal bundles \cite{bbn1}.
 For explicit descriptions we refer to \cite{furutasteer}, section 5,  \cite{boden}, section 4, or  \cite{biswasorbi}, section 2c.

 In \cite{nasatyrsteer2}, Nasatyr and Steer discuss Higgs bundles of rank $2$ in the orbifold setting, focusing on analytic aspects.  They define an orbifold Higgs bundle (or Higgs $V$-bundle) on an orbifold Riemann surface $X$ to be an orbifold bundle $E$ on $X$ together with an orbifold bundle morphism $\Phi:E\rightarrow E\otimes K$, where $K$ is the orbifold canonical bundle of $X$.  This definition is also used in \cite{nasatyrsteer1, bbn1}.  The correspondence with parabolic Higgs bundles is worked out in \cite{nasatyrsteer2}, section 5, as is the integrable system.  However, the parabolic Higgs bundles corresponding to the orbifold Higgs bundles they obtain are all (in our terminology) strongly parabolic (as one can observe by taking the residue of equation (5c) in \cite{nasatyrsteer2}).
 
 It would be interesting to discuss the matter of non-strongly parabolic Higgs bundles from an orbifold perspective and see if one could obtain the analogous results of our theorems \ref{tired}, \ref{endinsight}, \ref{groot}, \ref{oef}.  Presumably the analogue of all parabolic Higgs bundles (i.e. not necessarily strongly parabolic) would be given by looking at orbifold Higgs bundles with a Higgs field $E\rightarrow E\otimes L$, where $L$ is the orbifold line bundle obtained by twisting $K$ with (following the notation of \cite{nasatyrsteer2})the fractional divisor $\sum_i\frac{1}{\alpha_i}p_i$.
 An orbifold version of the work of Bottacin and Markman would then correspond to our results.

\appendix
\section{A Levi principal bundle over the moduli space of parabolic bundles}\label{append}

In this appendix we give a construction of a principal bundle  $\pfram \rightarrow \pvect$ with structure group $\mathcal{L}$, in the particular case of full flags at all marked points.  In the general case, one should think of the total space $\pfram$ of this principal bundle as a moduli space for  $\alpha$-stable parabolic bundles, together with isomorphisms of all consecutive quotients in the flags to a fixed vector space $$E_i(p)/E_{i+1}(p)\cong \mathbb{C}^{m_i}.$$ 
One can think of several approaches to this problem.
\\

One approach (for general flag types) one could take is suggested in \cite{hu}: start from a suitable moduli space of framed vector bundles (also known as bundles with level structure, in the case where the divisor over which one frames is reduced), as was for instance constructed in~\cite{huylehn}, generalizing earlier work by Seshadri~\cite{seshadri}.   The structure group of the vector bundles under consideration acts on this space by changing the framing, and we would like to take a GIT-style quotient by the unipotent radical of the parabolic subgroup.  For the Borel subgroup (leading to full flags) this is described in~\cite{hu}, a more general approach is given in~\cite{dorankirwan,frances}.  The various $\pvect$ would then be given by GIT quotients by the Levi group of this space, with the $\alpha$ occurring as the choice of a linearization.  As such the $\pvect$ are, for generic $\alpha$, geometric quotients for the Levi group actions, and therefore (with some mild extra conditions) principal bundles for the Levi group, see below.  
Notice that the action of $PGL(r)$ on a moduli space of framed bundles was discussed in~\cite{bisgommun}, where it was shown that the action linearizes and the GIT quotient is the moduli space of vector bundles.
\\

In order to keep the exposition from becoming too technical we shall use a construction already done in the literature, following~\cite{framparshea}.
Here a projective variety $\pmast$ was constructed directly, which was interpreted as a moduli space of framed parabolic sheaves\footnote{In~\cite{framparshea} the construction given further fixes the determinant for the sheaves, but the construction goes through without imposing this as well.}. The construction was inspired by a similar construction \cite{imp} in symplectic geometry through 
\emph{symplectic implosion}.  The connection between symplectic implosion and non-reductive GIT was discussed in \cite{preprintkirwan}. This variety $\mathcal{F}$ comes with a natural torus action, an action which linearizes on a relatively ample line bundle.  Using an earlier construction for the moduli space of parabolic bundles given by Bhosle~\cite{bhosle}, it is shown that at a linearization given by a character $\alpha$, the GIT quotient is the moduli space of parabolic vector bundles, $\pmast\git_{\alpha}T\cong \pvect$  (even in the case of partial flags, if one uses the $\alpha_i$ with the corresponding multiplicities).    If $\alpha$ is regular, i.e. we are looking at full flags, this is sufficient for us: the Luna slice theorem~\cite{luna}, see also \cite{GIT}, Appendix to Chapter 1, and \cite{shaves}, Corollary 4.2.13, now establishes that the $\alpha$-stable locus (the stability simply corresponds to the stability of the underlying bundle) $\pfram\subset\pmast$ is a principal bundle which is locally trivial in the \'etale topology. Furthermore, by a result of Serre~\cite{semchev}, as the structure group is a torus, it is even locally trivial in the Zariski topology.

\providecommand{\bysame}{\leavevmode\hbox to3em{\hrulefill}\thinspace}
\providecommand{\href}[2]{#2}

\end{document}